\documentclass{amsart}
\usepackage{amscd,amssymb,graphicx}
\author[DeCleene]{Chris DeCleene}
\address{
Chris DeCleene\\
University of Wisconsin-Eau Claire\\
Eau Claire, WI 54702-4004}
\email{cdecleene@gmail.com}
\author[Otto]{Carolyn Otto}
\address{
Carolyn Otto\\
Rice University\\
Houston, TX 77005-1827} \email{cotto@rice.edu}
\author[Penkava]{Michael Penkava}
\address{
Michael Penkava\\
University of Wisconsin-Eau Claire\\
Eau Claire, WI 54702-4004} \email{penkavmr@uwec.edu}
\author[Phillipson]{Mitch Phillipson}
\address{
Mitch Phillipson\\
University of Wisconsin-Eau Claire\\
Eau Claire, WI 54702-4004}
\email{phillima@uwec.edu}
\author[Steinbach]{Ryan Steinbach}
\address{
Ryan Steinbach\\
University of Wisconsin-Madison\\
Madison, WI 53706-1796}
\email{rsteinbach@wisc.edu}
\author[Weber]{Eric Weber}
\address{
Eric Weber\\
University of Wisconsin-Eau Claire\\
Eau Claire, WI 54702-4004}
\email{webered@uwec.edu}
\thanks{Research of these authors was partially supported by grants
from the National Science Foundation and the University of
Wisconsin-Eau Claire.}

\newtheorem{thm}{Theorem}[section]
\newtheorem{con}[thm]{Conjecture}

\theoremstyle{definition}

\def \ph{\varphi}

\def\GL{\mbox{\bf GL}}

\def \diag{\operatorname {diag}}

\def \ra{\rightarrow}

\def \hom{\mbox{\rm Hom}}
\def \ie{\hbox{\it i.e.}}

\def \tns{\otimes}

\def \mplus{+\cdots+}
\def \mcom{,\cdots,}
\def \k{\mbox{$\mathbb K$}}

\def \C{\mbox{$\mathbb C$}}
\def \Z{\mbox{$\mathbb Z$}}

\def\zt{\mbox{$\Z_2$}}

\def\ad{\operatorname{ad}}

\def\inv{^{-1}}

\def\im{\operatorname{Im}}

\def\m{\mbox{$\mathfrak m$}}

\def\coder{\operatorname{Coder}}
\def\ainf{\mbox{$A_\infty$}}

\def\and{\mbox{ \rm and }}
\def\T{\mathcal T}
\def\TA{\mbox{$\T(A)$}}
\def\TV{\mbox{$\T(V)$}}
\def\TW{\mbox{$\T(W)$}}

\def\s#1{(-1)^{#1}}
\DeclareMathOperator*{\invlim}{\underleftarrow{\rm lim}}

\def\qbar{\mbox{$\bar{\Q}$}}

\def\Q{\mbox{$\mathbb{Q}$}}

\def\pha#1#2{\ph^{#1}_{#2}}

\def\psa#1#2{\psi^{#1}_{#2}}

\def\inv{^{-1}}

\def\dinfty{\mbox{$d^\infty$}}

\def\P{\mathbb P}

\input{xy}
\xyoption{all}

\begin{document}
\setlength{\multlinegap}{0pt}
\title[The moduli space of $1|2$-dimensional complex associative algebras]
{Moduli space of $1|2$-dimensional complex associative algebras}%

\date{\today}
\begin{abstract}
In this paper, we study the moduli space of $1|2$-dimensional complex
associative algebras, which is also the moduli space of codifferentials
on the tensor coalgebra of a $2|1$-dimensional complex space.
We construct the moduli space by considering
extensions of lower dimensional algebras. We also construct miniversal
deformations of these algebras. This gives a complete description of
how the moduli space is glued together via jump deformations.
\end{abstract}
\maketitle

\nocite{ps1,ps2,mp1,pv1}
\section{Introduction}

The classification of associative algebras was instituted by Benjamin Peirce
in the 1870's \cite{pie}, who gave a partial classification of the complex
associative algebras of dimension up to 6, although in some sense, one can
deduce the complete classification from his results, with some additional
work. The classification method relied on the following remarkable fact:
\begin{thm}
Every finite dimensional algebra which is not nilpotent contains a nontrivial
idempotent element.
\end{thm}
A nilpotent algebra $A$ is one which satisfies $A^n=0$ for some $n$, while an
idempotent element $a$ satisfies $a^2=a$.  This observation of Peirce eventually
leads to two important theorems in the classification of finite dimensional associative
algebras.  Recall that an algebra is said to be simple if it has no nontrivial
proper ideals, and it is not the 1-dimensional nilpotent algebra over \k, given
by the trivial product.
\begin{thm}[Fundamental Theorem of Finite Dimensional Associative Algebras]
Suppose that $A$ is a finite dimensional algebra over a field \k. Then $A$ has
a maximal nilpotent ideal $N$, called its radical.  If $A$ is not nilpotent,
then $A/N$ is a semisimple
algebra, that is, a direct sum of simple algebras.
\end{thm}
Moreover, when $A/N$ satisfies a property called separability
over \k, then $A$ is a semidirect product of its radical and a semisimple algebra.
Over the complex numbers, every semisimple algebra is separable. To apply this theorem
to construct algebras by extension, one uses the following characterization of simple algebras.
\begin{thm}
[Wedderburn] If $A$ is a finite dimensional algebra over \k,
then $A$ is simple iff $A$ is isomorphic
to a tensor product $M\tns D$, where $M=\mathfrak{gl}(n,\k)$ and
$D$ is a division algebra over \k.
\end{thm}
In the nongraded case, a division algebra is a unital algebra
where every nonzero element has a multiplicative
inverse.
For \zt-graded associative algebras, the situation is a bit more complicated.
First, we need to
consider graded ideals, so that a \zt-graded algebra is simple
when it has no proper nontrivial
graded ideals. Secondly, the definition of a division algebra needs to be changed as well,
in order to generalize Wedderburn's theorem to the \zt-graded case. A \zt-graded division
algebra is a division algebra when every nonzero \emph{homogeneous} element is invertible.
With these changes, the Fundamental Theorem remains the same, except
that the radical is the maximal
graded nilpotent ideal, and Wedderburn's theorem is also true,
if we understand that by a matrix algebra,
we mean the general linear algebra of a \zt-graded vector space.

In this paper, we shall be concerned with the moduli space of associative algebras on a
\zt-graded algebra $A$ of dimension $2|1$, so that $A_0$ has
dimension 2 and $A_1$ has dimension 1.
However, we recall that this space coincides with the equivalence of odd
codifferentials of degree 2 on the
parity reversion $W=\Pi A$, which has dimension $1|2$, so in this paper,
we shall study codifferentials
on a space of dimension $1|2$, but the reader should keep in mind that
this corresponds to associative
algebras on a $2|1$-dimensional space.

The main goal of this paper is to give a complete description of
the moduli space of $2|1$-dimensional associative algebras,
including a computation of the miniversal deformation of every element.

\section{Construction of algebras by extensions}

In \cite{fp11}, the theory of extensions of an algebra $W$ by an algebra $M$ is
described in the language of codifferentials. Consider the diagram
$$
0\ra M\ra V\ra W\ra 0
$$
of associative \k-algebras, so that $V=M\oplus W$ as a \k-vector space, $M$ is an
ideal in the algebra $V$, and $W=V/M$ is the quotient algebra. Suppose that
$\delta\in C^2(W)$ and $\mu\in C^2(M)$ represent the algebra structures on
$W$ and $M$ respectively. We can view $\mu$ and $\delta$ as elements of $C^2(V)$.
Let $T^{k,l}$ be the subspace of $T^{k+l}(V)$ given recursively
by $T^{0,0}=\k$,
\begin{align*}
T^{k,l}&=M\tns T^{k-1,l}\oplus V\tns T^{k,l-1}.
\end{align*}
Let
$C^{k,l}=\hom(T^{k,l},M)\subseteq C^{k+l}(V)$.
If we denote the algebra structure on $V$ by $d$, we have
$$
d=\delta+\mu+\lambda+\psi,
$$
where $\lambda\in C^{1,1}$ and $\psi\in C^{0,2}$. Note that in this notation,
$\mu\in C^{2,0}$. Then the condition that $d$ is associative:  $[d,d]=0$ gives the
following relations:
\begin{align*}
[\delta,\lambda]+\tfrac 12[\lambda,\lambda]+[\mu,\psi]&=0,
\quad\text{The Maurer-Cartan equation}\\
[\mu,\lambda]&=0,\quad\text{The compatibility condition}\\
[\delta+\lambda,\psi]&=0,\quad\text{The cocycle condition}
\end{align*}
Since $\mu$ is an algebra structure, $[\mu,\mu]=0$, so if we define
$D_\mu$ by $D_\mu(\ph)=[\mu,\ph]$, then $D^2_\mu=0$.
Thus $D_\mu$ is a differential on $C(V)$.
Moreover $D_\mu:C^{k,l}\ra C^{k+1,l}$. Let
\begin{align*}
Z_\mu^{k,l}&=\ker(D_\mu:C^{k,l}\ra C^{k+1,l}),\quad\text{the $(k,l)$-cocycles}\\
B_\mu^{k,l}&=\im(D_\mu:C^{k-1,l}\ra C^{k,l}),\quad\text{the $(k,l)$-coboundaries}\\
H_\mu^{k,l}&=Z_\mu^{k,l}/B_\mu^{k,l},\quad\text{the $D_u$ $(k,l)$-cohomology}
\end{align*}

Then the compatibility condition means that $\lambda\in Z^{1,1}$.
If we define $D_{\delta+\lambda}(\ph)=[\delta+\lambda,\ph]$, then it is not
true that $D^2_{\delta+\lambda}=0$, but
$D_{\delta+\lambda}D_\mu=-D_{\mu}D_{\delta+\lambda}$, so that $D_{\delta+\lambda}$ descends
to a map $D_{\delta+\lambda}:H^{k,l}_\mu\ra H^{k,l+1}_\mu$, whose square is zero, giving
rise to the $D_{\delta+\lambda}$-cohomology $H^{k,l}_{\mu,\delta+\lambda}$.
If the pair $(\lambda,\psi)$ give rise to a codifferential $d$, and $(\lambda,\psi')$
give rise to another codifferential $d'$, then if we express $\psi'=\psi+\tau$, it is
easy to see that $[\mu,\tau]=0$, and $[\delta+\lambda,\tau]=0$, so that the image $\bar\tau$
of $\tau$ in $H^{0,2}_\mu$ is a $D_{\delta+\lambda}$-cocycle, and thus $\tau$ determines
an element $\{\bar\tau\}\in H^{0,2}_{\mu,\delta+\lambda}$.

If $\beta\in C^{0,1}$, then $g=\exp(\beta):\TV\ra\TV$ is given by
$g(m,w)=(m+\beta(w),w)$. Furthermore $g^*=\exp(-\ad_{\beta}):C(V)\ra C(V)$ satisfies
$g^*(d)=d'$, where $d'=\delta+\mu+\lambda'+\psi'$ with
\begin{align*}
\lambda'&=\lambda+[\mu,\beta]\\
\psi'&=\psi+[\delta+\lambda+\tfrac12[\mu,\beta],\beta],
\end{align*}
In this case, we say that $d$ and $d'$ are equivalent extensions in the restricted sense.
Such equivalent extensions are also equivalent as codifferentials on $\TV$.
Note that $\lambda$
and $\lambda'$ differ by a $D_\mu$-coboundary, so $\bar\lambda=\bar\lambda'$ in
$H^{1,1}_\mu$. If $\lambda$ satisfies the MC-equation for some $\psi$, then
any element $\lambda'$ in $\bar\lambda$ also gives a solution of the MC equation,
for the $\psi'$ given above. The cohomology classes of those $\lambda$ for which
a solution of the MC equation exists determine distinct restricted equivalence classes
of extensions.

Let $G_{M,W}=\GL(M)\times\GL(W)\subseteq\GL(V)$. If $g\in G_{M,W}$ then $g^*:C^{k,l}\ra
C^{k,l}$, and $g^*:C^k(W)\ra C^k(W)$, so $\delta'=g^*(\delta)$ and $\mu'=g^*(\mu)$
are codifferentials on $\T(M)$ and $\TW$ respectively.
The group $G_{\delta,\mu}$ is the
subgroup of $G_{M,W}$ consisting of those elements $g$ such that $g^*(\delta)=\delta$
and $g^*(\mu)=\mu$. Then $G_{\delta,\mu}$ acts on
the restricted equivalence classes of extensions, giving the equivalence classes
of general extensions. Also, $G_{\delta,\mu}$
acts on $H^{k,l}_\mu$, and induces an action on the classes $\bar\lambda$ of $\lambda$
giving a solution $(\lambda,\psi)$ to the MC equation.

Next, consider the group $G_{\delta,\mu,\lambda}$ consisting
of the automorphisms $h$ of $V$ of the form $h=g\exp(\beta)$, where
$g\in G_{\delta,\mu}$, $\beta\in C^{0,1}$ and $\lambda=g^*(\lambda)+[\mu,\beta]$.
If $d=\delta+\mu+\lambda+\psi+\tau$, then $h^*(d)=\delta+\mu+\lambda+\psi+\tau'$ where
\begin{equation*}
\tau'=g^*(\psi)-\psi+[\delta+\lambda-\tfrac12[\mu,\beta],\beta]+g^*(\tau).
\end{equation*}
Moreover, the group $G_{\delta,\mu,\lambda}$ induces an action on $H^{0,2}_{\mu,\delta+\lambda}$
given by $\{\bar\tau\}\ra\{\overline{\tau'}\}$. In fact, $\{\overline{g^*(\tau)}\}$ is well
defined as well, and depends only on $\{\bar\tau\}$.

The general group of equivalences of extensions of the algebra structure $\delta$ on $W$
by the algebra structure $\mu$ on $M$ is given by the group of automorphisms of $V$ of
the form $h=\exp(\beta)g$, where $\beta\in C^{0,1}$ and $g\in G_{\delta,\mu}$. We
have the following classification of such extensions up to equivalence.
\begin{thm}
The equivalence classes of extensions of $\delta$ on $W$ by $\mu$ on $M$ is classified
by the following:
\begin{enumerate}
\item Equivalence classes of $\bar\lambda\in H^{1,1}_\mu$ which satisfy the MC equation
\begin{equation*}
[\delta,\lambda]+\tfrac12[\lambda,\lambda]+[\mu,\psi]=0
\end{equation*}
for some $\psi\in C^{0,2}$, under the action of the group $G_{\delta,\mu}$.
\item Equivalence classes of $\{\bar\tau\}\in H^{0,2}_{\mu,\delta+\lambda}$ under the
action of the group $G_{\delta,\mu,\lambda}$.
\end{enumerate}
\end{thm}
Equivalent extensions will give rise to equivalent codifferentials on $V$, but it may
happen that two codifferentials arising from nonequivalent extensions are equivalent.
This is because the group of equivalences of extensions is the group of invertible
block upper
triangular matrices on the space $V=M\oplus W$, whereas the equivalence
classes of codifferentials on $V$ are given by the group of all invertible
matrices, which is larger.

The fundamental theorem of finite dimensional algebras allows us to restrict our
consideration of extensions to two cases. First, we can consider those extensions
where $\delta$ is a semisimple algebra structure on $W$, and $\mu$ is a nilpotent
algebra structure on $M$. In this case, because we are working over $\C$, we can
also assume that $\psi=\tau=0$. Thus the classification of the extension reduces
to considering equivalence classes of $\lambda$.

Secondly, we can consider extensions
of the trivial algebra structure $\delta=0$ on a 1-dimensional space $W$ by
a nilpotent algebra $\mu$. This
is because a nilpotent algebra has a codimension 1 ideal $M$, and the restriction
of the algebra structure to $M$ is nilpotent. However, in this case, we cannot assume
that $\psi$ or $\tau$ vanish,
so we need to use the classification theorem above to determine the
equivalence classes of extensions. In many cases, in solving the MC equation for
a particular $\lambda$, if there is any $\phi$ yielding a solution, then $\psi=0$
also gives a solution, so the action of $G_{\delta,\mu,\lambda}$ on $H^{0,2}_\mu$
takes on a simpler form than the general action we described above. In fact, if in
addition to $\psi=0$ providing a solution to the MC equation, any element $h=g\exp(\beta)$
satisfies $[\mu,\beta]=0$, then the action of $h$ on $H^{0,2}_{\delta,\mu,\lambda}$ is
just the action $g^*(\{\bar\tau\})=\{\overline{g^*(\tau)}\}$, which is easy to calculate
in practice.

\section{Associative algebra structures on a $1|2$ vector space}
Let $A$ be a $1|2$-dimensional vector space, and $V=\Pi A$ be the
parity reversion of $A$, so that $V$ is $2|1$-dimensional.
Let $\{v_1,v_2,v_3\}$ be a basis of $V$ with $v_1, v_2$ even
elements and $v_3$ an odd element, and let $d$ be a codifferential
on $V$ representing an associative algebra structure on $A$.

By results in
\cite{bppw1}, there are only two
\zt-graded division algebras, the complex numbers, and a certain $1|1$-dimensional
algebra. As a consequence, there are no $1|2$-dimensional simple algebras,
so we can express $V$ as an extension of an algebra structure $\delta$ on $W$
by an algebra structure $\mu$ on $M$, where $V=M\oplus W$, and $M$ is an ideal
in $V$.

By the fundamental theorem of finite dimensional associative algebras,
we can assume that $\mu$ is a nilpotent algebra structure on $M$. Moreover,
$\delta$ is a semisimple algebra structure on $W$, unless
$d$ is a nilpotent algebra structure.

Since every nilpotent algebra has a codimension 1 ideal, if $d$ is nilpotent,
we can assume that $W$ is 1-dimensional (either even or odd), and that $\delta=0$.
The only semisimple algebras we need to consider are the simple $1|1$-dimensional
algebra, and the simple $0|1$-dimensional algebra. Moreover, when considering extensions
of a semisimple algebra, the ``cocycle`` $\psi$ can be taken to be zero,
because we are considering extensions over $\C$,
for which every semisimple algebra is separable.

Now,suppose that $W=\langle v_{w(1)}\mcom v_{w(p)}\rangle$ where the first $s$ vectors are even
and the other $p-s$ elements are odd, and that
$M=\langle v_{m(1)}\mcom v_{m(q)}\rangle$, where the first $t$ elements are even and
the other $q-t$ elements are odd. (This conforms to the principle that in a \zt-graded space,
a basis should be listed with the even elements first.) Then a formula for an arbitrary
$\lambda\in C^{1,1}$ is of the form
\def\LE{LE}
\begin{align*}
\lambda&=\sum_{|v_{w(k)}|=1}\psa{w(k)m(j)}{m(i)}(\LE_k)^i_j+\psa{m(j)w(k)}{m(i)}(RE_k)^i_j
\\&+\sum_{|v_{w(k)}|=0} \psa{w(k)m(j)}{m(i)}(LO_k)^i_j+\psa{m(j)w(k)}{m(i)}(RO_k)^i_j,
\end{align*}
where $LE_k$ and $RE_k$ are matrices of even maps $M\ra M$, and $LO_k$ and $RO_k$ are matrices
of odd maps $M\ra M$. For simplicity, we shall denote $(LE_k)^i_J$ as $LE_{kj}^i$ and
similarly for the components of the other matrices. Let $L_k=LE_k+LO_k$ and $R_k=RE_k+RO_k$.
Then
\begin{align*}
\tfrac12[\lambda,\lambda]&=\,\psa{w(k)w(l)m(j)}{m(i)}((LO_k-LE_k)L_l))^i_j+
\psa{m(j)w(k)w(l)}{m(i)}(R_lR_k)^i_j
\\&\,+\psa{w(k)m(j)w(l)}{w(i)}(R_lL_k+(LO_k-LE_k)R_l).
\end{align*}
It is important to note that the formula above is given in terms of matrix multiplication.
This is significant from the computational view as we shall illustrate below. It is interesting
to note that matrices in $G_{M,W}$, which are block diagonal maps $\diag(G_1,G_2)$, where
$G_1\in\GL(M)$ and $G_2\in\GL(W)$, act on $\lambda$ in a manner which can be described in
terms of the matrices above. First, $G_1$ acts by conjugating all the matrices $L_k$ and
$R_k$ simultaneously. Secondly, the matrix $G_2$ acts on the $k$ indices. We shall say
more about these actions later.

Let us give one concrete application of the remarks above. Suppose that $W$ is completely odd,
and $M$ is $r|s$-dimensional. Then we can express $M=\langle v_1\mcom v_{r+s}\rangle$ and
$W=\langle v_{r+s+1}\mcom v_{r+s+n}\rangle$. In this case $L_k=LE_k$ and $R_k=RE_k$.
We can express $R_k=\diag(T_k,B_k)$, where $T_k:M_0\ra M_0$ and $B_k:M_1\ra M_1$ and
$M=M_0\oplus M_1$ represents the decomposition of $M$ into its even and odd parts.
Then we have
\begin{equation*}
\tfrac12[\lambda,\lambda]=-\psa{klj}i(L_kL_l)^i_j+\psa{kjl}i(R_lL_k-L_kR_l)^i_j
+\psa{jkl}i(R_lR_k)^i_j.
\end{equation*}
Let $\delta=\sum_{k=r+s+1}^{r+s+n}\psa{kk}k$ be  the codifferential semisimple algebra $\C^n$.
Then
\begin{equation*}
[\delta,\lambda]=\psa{kkj}i(L_{k})^i_j+\s{j}\psa{jkk}i(R_{k})^i_j.
\end{equation*}
When considering an extension of the algebra $\C^n$ by an algebra structure on $M$, we
can assume that the cocycle ``$\psi$`` vanishes, so the MC equation is just
$[\delta,\lambda]+\tfrac12[\lambda,\lambda]=0$, which is equivalent to the following:
\begin{center}
\begin{align*}
L_k=L_k^2,\quad &L_kL_l=0,\text{ if $k\ne l$}\\
L_k&R_l=R_lL_k\\ T_k=-T_k^2,\quad B_k=&B_k^2,\quad R_kR_l=0,\text{ if $k\ne l$}.
\end{align*}
\end{center}
As a consequence, the matrices above give a commuting set of diagonalizable matrices,
so they can be simultaneously diagonalized. Moreover, $L_k$ and $B_k$ have only
0 and 1 as possible eigenvalues and $T_k$ has only 0 and $-1$ as possible eigenvalues.
If we consider $\mu=0$ as the algebra structure on $M$, then the elements
$\diag(G_1,G_2)\in G_{\delta,\mu}$ are given by an arbitrary matrix $G_1$ and $G_2$ is
just a permutation matrix.

Thus we can apply an element $G\in G_{\delta,\mu}$ to
$\lambda$ to put it in the form where all the matrices are diagonal, and are ordered
in such a manner that the nonzero $L_k$ matrices appear first. Moreover, since the
$L_k$ matrices are orthogonal to each other, there are at most $m=r+s$ nonzero $L$
matrices.

Similar considerations apply to the $R$ matrices, so that in total, there
can be no more than $2m$ pairs $(L_k,R_k)$, such that at least one matrix does not
vanish.  Therefore, when $n\ge 2m$, the number of distinct equivalence classes of
extensions of $\C^n$ by a trivial algebra structure on $M$ is exactly equal to $2m$.

We say that this is the stable situation. The number of extensions is independent
of $n$ as long as it is at least $2m$. Moreover, the cohomology and deformation
theory also becomes stable, in a natural way. When $\mu\ne0$, the situation is a
bit more complicated, but there is also an $n$ beyond which the situation becomes
stable.

We now give a construction of the elements in the moduli space of $2|1$-dimensional
codifferentials. Table \ref{coho21 table} below gives the cohomology of the 21
nonequivalent codifferentials.

\begin{table}[h!]
\begin{center}
\begin{tabular}{lccccc}
Codifferential&$H^0$&$H^1$&$H^2$&$H^3$&$H^4$\\ \hline \\
$d_1=\psa{23}2-\psa{32}2+\psa{22}3-\psa{33}3$&$1|1$&$1|0$&$1|0$&$1|0$&$1|0$\\
$d_2=\psa{33}3+\psa{31}1+\psa{32}2$&$0|0$&$3|0$&$0|0$&$0|0$&$0|0$\\
$d_3=\psa{33}3-\psa{13}1-\psa{23}2$&$0|0$&$3|0$&$0|0$&$0|0$&$0|0$\\
$d_4=\psa{33}3+\psa{31}1-\psa{23}2$&$0|0$&$1|0$&$0|0$&$1|0$&$0|0$\\
$d_5=\psa{33}3-\psa{13}1+\psa{32}2-\psa{23}2$&$1|0$&$1|0$&$1|0$&$1|0$&$1|0$\\
$d_6=\psa{33}3+\psa{31}1+\psa{32}2-\psa{23}2$&$1|0$&$1|0$&$1|0$&$1|0$&$1|0$\\\hline
$d_7=\psa{33}3+\psa{32}2$&$1|0$&$1|0$&$2|0$&$2|0$&$2|0$\\
$d_8=\psa{33}3-\psa{23}2$&$1|0$&$1|0$&$2|0$&$2|0$&$2|0$\\
$d_9=\psa{33}3+\psa{31}1-\psa{13}1+\psa{32}2-\psa{23}2$&$3|0$&$4|0$&$6|0$&$12|0$&$24|0$\\
$d_{10}=\psa{33}3$&$3|0$&$4|0$&$8|0$&$16|0$&$32|0$\\
$d_{11}=\psa{33}3+\psa{32}2-\psa{23}2$&$2|1$&$2|1$&$2|1$&$2|1$&$2|1$\\
$d_{12}=\psa{22}3+\psa{23}1-\psa{32}1$&$1|1$&$2|0$&$1|1$&$2|0$&$1|1$\\\hline
$d_{13}(p:q)=\psa{22}3+\psa{21}3p+\psa{12}3q$&$0|1$&$2|0$&$2|1$&$3|0$&$4|0$\\
$d_{13}(1:1)=\psa{22}3+\psa{21}3+\psa{12}3$&$0|1$&$2|0$&$2|1$&$5|0$&$4|2$\\
$d_{13}(1:-1)=\psa{22}3+\psa{21}3-\psa{12}3$&$1|1$&$2|1$&$3|1$&$4|1$&$5|1$\\
$d_{13}(1:0)=\psa{22}3+\psa{21}3$&$1|0$&$2|0$&$4|1$&$6|2$&$8|3$\\
$d_{13}(0:0)=\psa{22}3$&$1|1$&$3|1$&$5|4$&$10|7$&$18|14$\\\hline
$d_{14}=\psa{21}3-\psa{12}3$&$2|1$&$4|2$&$5|4$&$8|4$&$10|5$\\
$d_{15}(p:q)=\psa{23}1p+\psa{32}1q$&$1|0$&$2|0$&$1|2$&$2|1$&$2|2$\\
$d_{15}(1:1)=\psa{23}1+\psa{32}1$&$1|0$&$2|1$&$2|2$&$4|2$&$3|4$\\
$d_{15}(1:0)=\psa{23}1$&$1|0$&$2|0$&$2|3$&$5|3$&$5|6$\\
$d_{15}(0:1)=\psa{32}1$&$1|0$&$2|0$&$2|3$&$5|3$&$5|6$\\
$d_{15}(1:-1)=\psa{23}1-\psa{32}1$&$2|1$&$3|2$&$4|3$&$5|4$&$6|5$\\
\\ \hline
\end{tabular}
\end{center}
\label{coho21 table}
\caption{Cohomology of the 15 families of codifferentials on a $2|1$-dimensional space}
\end{table}

\section{Extensions where $W$ is $1|1$-dimensional and $M$ is $1|0$-dimensional}
Let $W=\langle v_2,v_3\rangle$ and $M=\langle v_1\rangle$.
The unique $1|1$-dimensional simple algebra is given by the codifferential
$\delta=\psa{23}2-\psa{32}2+\psa{22}3-\psa{33}3$. The only algebra structure
on $M$ is the trivial algebra $\mu=0$.  The generic lambda is of the form
$\lambda=\psa{31}1LE_{21}^1+\psa{13}1RE_{21}^1$. However,
\begin{align*}
[\delta,\lambda]&=\pha{221}1LE_{21}^1-\pha{122}1RE_{21}^1
-\pha{331}1LE_{21}^1-\pha{133}1RE_{21}^1\\
\tfrac12[\lambda,\lambda]&=-\pha{331}1(LE_{21}^1)^2+\pha{133}1(RE_{21}^1)^2
\end{align*}
so the MC equation forces $\lambda=0$. Therefore, the unique extension of $\delta$
is the direct sum of $\delta$ and the trivial 1-dimensional algebra, which
is the codifferential $d_1$.

\section{Extensions where $M$ is $2|0$-dimensional and $W$ is $0|1$-dimensional}
Let $M=\langle v_1,v_2\rangle$ be $2|0$-dimensional and $W=\langle v_3\rangle$
be $0|1$-dimensional.
The group $G_{M,W}$ consists of diagonal matrices
$G=\diag(r,s,t)$, with $rst\ne0$. A
generic element of $C^{1,1}$ is of the form
\begin{align*}
\lambda&=
\psa{31}1aL_{1}^1+\psa{31}2L_{2}^1+\psa{32}1L_2^1+\psa{32}2L_2^2
\psa{13}1aR_{1}^1+\psa{13}2R_{2}^1+\psa{23}1R_2^1+\psa{23}2R_2^2
\end{align*}
corresponding to the matrices
\begin{equation*}
L_{{1}}= \left[ \begin {matrix}
L_1^1&L_2^1\\\noalign{\medskip}
L_1^2&L_2^2
\end{matrix} \right] ,
R_{{1}}=
\left[ \begin {matrix}
R_1^1&R_2^1\\\noalign{\medskip}
R_1^2&R_2^2
\end {matrix} \right].
\end{equation*}

\subsubsection{Extensions of the simple $0|1$-dimensional algebra}
In this case $\delta=\psa{33}3$. From the MC-equation, we obtain that
\begin{equation*}
L^2=L,\quad RL=LR,\quad R^2=-R,
\end{equation*}
so that $L$ is a diagonalizable matrix with eigenvalues 1 and 0, and
$R$ is diagonalizable with eigenvalues $-1$ and 0. Thus $R$ and $L$ can
be simultaneously diagonalized, which means that we have the following
solutions. If $L=I$ or $L=0$, then $R$ can be taken to be either $-I$, $\diag(-1,0)$ or $0$.
Otherwise $L=\diag(1,0)$ and $R$ is either $I$, $\diag(-1,0)$, $\diag(0,-1)$, or 0.
These 10 possibilities correspond to the codifferentials $d_2$,\dots $d_{11}$.

\subsubsection{Extensions of the trivial $0|1$-dimensional algebra}.
In this case, $\delta=0$ and $\mu=0$. Then $G_{\delta,\mu}=G_{M,W}$.
The MC-equation yields $L^2=0$, $LR=RL$ and $R^2=0$. The eigenvalues
of the matrices $R$ and $L$ are only 0, and since these matrices commute,
they can put in simultaneous upper triangular form. Thus we can express
$\lambda=\psa{32}2p+\psa{23}2q$, which gives a family of nonequivalent
codifferentials, $d_{15}(p:q)$. This family is
parameterized projectively by elements $(p:q)\in\C\P^1$.

\section{Extensions where $M$ is $1|1$-dimensional and $W$ is $1|0$-dimensional}
We have $M=\langle v_1,v_3\rangle$ and $W=\langle v_2\rangle$.
The group $G_{M,W}$ consists of matrices of the form
$G=\diag(r,s,t)$, where $rst\ne0$. The only codifferential on $W$ is $\delta=0$.
A generic element in $C^{1,1}$ is
\begin{equation*}
\lambda=\psa{23}1LO^1_2+\psa{21}3LO^2_1
+\psa{32}1RO^1_2+\psa{12}3RO^2_1,
\end{equation*}
corresponding to the matrices
\begin{equation*}
L= \left[ \begin {array}{cc} 0&LO^1_2
\\\noalign{\medskip}LO^2_1&0\end {array} \right] ,R=
\left[ \begin {array}{cc} 0&RO^1_2\\\noalign{\medskip}RO^2_1&0
\end {array} \right] ]
\end{equation*}
Moreover a generic element in $C^{0,2}$ is of the form $\tau=\psa{22}3c$
and a generic element of $C^{0,1}$ is of the form $\beta=\pha21b$.
There are two
possibilities for $\mu$, the trivial codifferential and
$\mu=\psa{11}2$.

In this case, it is more convenient to consider the trivial case $\mu=0$ first.
Then the MC equation yields that either  $LO^2_1=RO^2_1=0$ or $LO^1_2=RO^1_2=0$.

Let us consider the first case. Let $LO^1_2=p$ and $RO^1_2=q$.
Then $\lambda=\psa{23}1p+\psa{32}1q$. Now $[\delta+\lambda,\tau]=\pha{222}1(p+q)c$,
so unless $p=-q$, we can assume that $\tau=0$. Therefore, unless $p=-q$, we
have $d=\psa{23}1p+\psa{32}1q$, which is just the codifferential $d_{15}(p: q)$ we
already encountered. We also note that the action of $G_{\delta,\mu}$ on $\lambda$ is
to multiply it by a nonzero constant, so that $\lambda$ can be parameterized projectively.
Thus if $p=-q$, we can assume that either $p=1$ and $q=-1$ or $p=q=0$. In the
either case, the action of $G_{\delta,\mu,\lambda}$ on $\tau$ multiplies it by a constant,
so we can assume that $\tau=\psa{22}3$, since the case $\tau=0$ is generic.
This gives the codifferential $d_{12}$ when $\lambda=\psa{23}1-\psa{32}1$, and
the codifferential $d_{13}(0:0)$ when $\lambda=0$.

In the second case, we can take $\lambda=\psa{21}3p+\psa{12}3q$. Then
$G_{\delta,\mu,\lambda}$ consists of those transformations of the form
$h=g\exp(\beta)$ where $\beta=\pha{2}1b$ and g is given by a matrix $G=\diag(r,s,rs)$.
Note that $[\delta+\lambda,\beta]=\psa{22}3(p+q)b$, which means that all values of $c$
produce an equivalent codifferential unless $p=-q$.

It might seem that this would
dictate choosing $c=0$ in the generic case, but this actually is not what deformation
theory requires. The reason is that in the special cases where $p=-q$, for all values
of $c$ except $c=0$ we obtain an equivalent codifferential, so the codifferential
arising from taking $c=0$ jumps to the one arising from taking $c=1$. Therefore,
the case $c=1$ is generic.

In the generic case, we take $c=1$ and we obtain the codifferential $d_{13}(p:q)$,
which is parameterized projectively by $(p:q)\in\C\P^1$. In the case where $p=1$ and
$q=-1$ and $\tau=0$, we obtain the codifferential $d_{14}$, and when $p=q=0$ and
$\tau=0$, we obtain the zero codifferential.

Now, consider the nontrivial case $\mu=\psa{11}3$.
Then $G_{\delta,\mu}$ is given by matrices of the form $\diag(r,s,s^2)$ such
that $rs\ne0$.
The compatibility condition $[\mu,\lambda]=0$ yields
$LO^2_1=0$ and $RO^2_1=0$, and $[\mu,\beta]=\psa{21}3b+\psa{12}3b$, so
we can assume that $RO^1_2=0$ as well. Taking in to account the action
of the group $G_{\delta,\mu}$, we can reduce to the cases $\lambda=\psa{21}3$ or
$\lambda=0$.

If $\lambda=\psa{21}3$, then the action of $G_{\delta,\mu,\lambda}$ on $\tau$ leaves
it unchanged, so we have to consider all values of $c$. It turns out that
the codifferentials we obtain are equivalent to
$d_{13}(-1+\sqrt{1-4c}:1+\sqrt{1-4c})$. This complicated formula is why it was more
convenient to study the case $\mu=0$ first.

If $\lambda=0$, then the action of $G_{\delta,\mu,\lambda}$ on $\tau$ multiplies it
by a nonzero constant, so we can consider the case $c_1=0$, which gives a codifferential
equivalent to $d_{13}(0:0)$, or $c_1=1$, which gives the codifferential $d_{13}(1:1)$.

Thus we have completed the classification of the elements in the moduli
space of associative algebra structures on a space $V$ of dimension $1|2$,
or, in other words, the codifferentials on a $2|1$-dimensional space.

\section{Hochschild Cohomology and Deformations}
\emph{Hochschild cohomology} was introduced in \cite{hoch}, and used to
classify infinitesimal deformations of associative algebras. Suppose that
\begin{equation*}
m_t=m+t\ph,
\end{equation*}
is an infinitesimal deformation of $m$.  By this we mean that the structure
$m_t$ is associative up to first order. From an algebraic point of view, this
means that we assume that $t^2=0$, and then check whether associativity holds.
It is not difficult to show that is equivalent to the following.
\begin{equation*}
a\ph(b,c)-\ph(ab,c)+\ph(a,bc)-\ph(a,b)c=0,
\end{equation*}
where, for simplicity, we denote $m(a,b)=ab$. Moreover, if we let
\begin{equation*}
g_t=I+t\lambda
\end{equation*} be an infinitesimal automorphism of $A$, where
$\lambda\in\hom(A,A)$, then it is easily checked that
\begin{equation*}
g_t^*(m)(a,b)=ab+t(a\lambda(b)-\lambda(ab)+\lambda(a)b).
\end{equation*}
 This naturally leads
to a definition of the Hochschild coboundary operator $D$ on $\hom(\TA,A)$ by
\begin{align*}
D(\ph)(a_0\mcom a_n)=&a_0\ph(a_1\mcom a_n)+\s{n+1}\ph(a_0\mcom a_{n-1})a_n\\
&+\sum_{i=0}^{n-1}\s{i+1}\ph(a_0\mcom a_{i-1},a_ia_{i+1},a_{i+2}\mcom a_n)
.
\end{align*}
If we set $C^n(A)=\hom(A^n,A)$, then $D:C^n(A)\ra C^{n+1}(A)$. One obtains the
following classification theorem for infinitesimal deformations.
\begin{thm} The equivalence classes of infinitesimal
deformations $m_t$ of an associative algebra structure $m$ under the action of the
group of infinitesimal automorphisms on the set of infinitesimal deformations
are classified by the Hochschild cohomology group
\begin{equation*}
H^2(m)=\ker(D:C^2(A)\ra C^3(A))/\im(D:C^1(A)\ra C^2(A)).
\end{equation*}
\end{thm}
When $A$ is \zt-graded, the  only modifications that are necessary are that
$\ph$ and $\lambda$ are required to be even maps, so we obtain that the
classification is given by $H^2_e(A)$, the even part of the Hochschild
cohomology.

We wish to transform this classical viewpoint into the more modern viewpoint of
associative algebras as being given by codifferentials on a certain coalgebra.
To do this, we first introduce the \emph{parity reversion} $\Pi A$ of a
\zt-graded vector space $A$. If $A=A_e\oplus A_o$ is the decomposition of $A$
into its even and odd parts, then $W=\Pi A$ is the \zt-graded vector space
given by $W_e=A_o$ and $W_o=A_e$. In other words, $W$ is just the space $A$
with the parity of elements reversed.

Denote the tensor (co)-algebra of $W$ by $\TW=\bigoplus_{k=0}^\infty W^k$,
where $W^k$ is the $k$-th tensor power of $W$ and $W^0=\k$. For brevity, the
element in $W^k$ given by the tensor product of the elements $v_i$ in $W$ will
be denoted by $v_1\cdots v_k$. The coalgebra structure on $\TW$ is  given by
\begin{equation*}
\Delta(v_1\cdots v_n)=\sum_{i=0}^n v_1\cdots v_i\tns v_{i+1}\cdots v_n.
\end{equation*}
Define $d:W^2\ra W$ by $d=\pi\circ m\circ (\pi\inv\tns\pi\inv)$, where
$\pi:A\ra W$ is the identity map, which is odd, because it reverses the parity
of elements. Note that $d$ is an odd map. The space $C(W)=\hom(\TW,W)$ is
naturally identifiable with the space of coderivations of $\TW$.  In fact, if
$\ph\in C^k(W)=\hom(W^k,W)$, then $\ph$ is extended to a coderivation of $\TW$
by
\begin{equation*}
\ph(v_1\cdots v_n)=
\sum_{i=0}^{n-k}\s{(v_1\mplus v_i)\ph}v_1\cdots
 v_i\ph(v_{i+1}\cdots v_{i+k})v_{i+k+1}\cdots v_n.
\end{equation*}

The space of coderivations of $\TW$ is equipped with a \zt-graded Lie algebra
structure given by
\begin{equation*}
[\ph,\psi]=\ph\circ\psi-\s{\ph\psi}\psi\circ\ph.
\end{equation*}
The reason that it is more convenient to work with the structure $d$ on $W$
rather than $m$ on $A$ is that the condition of associativity for $m$
translates into the codifferential property $[d,d]=0$.  Moreover, the
Hochschild coboundary operation translates into the coboundary operator $D$ on
$C(W)$, given by
\begin{equation*}
D(\ph)=[d,\ph].
\end{equation*}
This point of view on Hochschild cohomology first appeared in \cite{sta4}.  The
fact that the space of Hochschild cochains is equipped with a graded Lie
algebra structure was noticed much earlier \cite{gers,gers1,gers2,gers3,gers4}.

For notational purposes, we introduce a basis of $C^n(W)$ as follows.  Suppose
that $W=\langle v_1\mcom v_m\rangle$. Then if $I=(i_1\mcom i_n)$ is a
\emph{multi-index}, where $1\le i_k\le m$, denote $v_I=v_{i_1}\cdots v_{i_n}$.
Define $\ph^{I}_i\in C^n(W)$ by
\begin{equation*}
\ph^I_i(v_J)=\delta^I_Jv_i,
\end{equation*}
where $\delta^I_J$ is the Kronecker delta symbol. In order to emphasize the
parity of the element, we will denote $\ph^I_i$ by $\psi^I_i$ when it is an odd
coderivation.

For a multi-index $I=(i_1\mcom i_k)$, denote its \emph{length}  by $\ell(I)=k$.  If
$K$ and $L$ are multi-indices, then denote $KL=(k_1\mcom k_{\ell(K)},l_l\mcom
l_{\ell(L)})$.  Then
\begin{align*}
(\ph^I_i\circ\ph^J_j)(v_K)&=
\sum_{K_1K_2K_3=K}\s{v_{K_1}\ph^J_j} \ph^I_i(v_{K_1},\ph^J_j(v_{K_2}), v_{K_3})
\\&=
\sum_{K_1K_2K_3=K}\s{v_{K_1}\ph^J_j}\delta^I_{K_1jK_3}\delta^J_{K_2}v_i,
\end{align*}
from which it follows that
\begin{equation}\label{braform}
\ph^I_i\circ\ph^J_j=\sum_{k=1}^{\ell(I)}\s{(v_{i_1}\mplus v_{i_{k-1}})\ph^J_j}
\delta^k_j
\ph^{(I,J,k)}_i,
\end{equation}
where $(I,J,k)$ is given by inserting $J$ into $I$ in place of the $k$-th
element of $I$; \ie, $(I,J,k)=(i_1\mcom i_{k-1},j_1\mcom j_{\ell(J)},i_{k+1}\mcom
i_{\ell(I)})$.

Let us recast the notion of an infinitesimal deformation in terms of the
language of coderivations.  We say that
\begin{equation*}
d_t=d+t\psi
\end{equation*}
is a deformation of the codifferential $d$ precisely when $[d_t,d_t]=0 \mod t^2$.
This condition immediately reduces to the cocycle condition $D(\psi)=0$.  Note
that we require $d_t$ to be odd, so that $\psi$ must be an odd coderivation.
One can introduce a more general idea of parameters, allowing both even and odd
parameters, in which case even coderivations play an equal role, but we will
not adopt that point of view in this paper.

For associative algebras, we require that $d$ and $\psi$ lie in $\hom(W^2,W)$.
This notion naturally generalizes to considering $d$ simply to be an arbitrary
odd codifferential, in which case we would obtain an \ainf\ algebra, a natural
generalization of an associative algebra.

More generally, we need the notion of a versal deformation, in order to
understand how the moduli space is glued together. To explain versal deformations we introduce
the notion of  a deformation with a local base.
A local base $A$ is a \zt-graded commutative, unital
$\k$-algebra with an augmentation $\epsilon:A\ra\k$,
whose kernel $\m$ is the unique maximal ideal in $A$,
so that $A$ is a local ring. It follows that $A$ has a unique decomposition
$A=\k\oplus\m$ and $\epsilon$ is
just the projection onto the first factor. Let $W_A=W\tns A$ equipped
with the usual structure of a right $A$-module.
Let $T_A(W_A)$ be the tensor algebra of $W_A$ over $A$,
that is $T_A(W_A)=\bigoplus_{k=0}^\infty T^k_A(W_A)$ where $T^0_A(W_A)=A$ and
$T^{k+1}_A(W_A)=T^k(W_A)_A\tns_A W_A$.
It is a standard fact that $T^k_A(W_A)=T^k(W)\tns A$ in a natural manner,
and thus $T_A(W_A)=T(W)\tns A$.

Any $A$-linear map $f:T_A(W)\ra T_A(W)$ is induced by its restriction to
$T(W)\tns \k=T(W)$ so we can view an $A$-linear coderivation
$\delta_A$ on $T_A(W_A)$ as a map $\delta_A:T(W)\ra T(W)\tns A$.
A morphism $f:A\ra B$ induces a map  $$f_*:\coder_A(T_A(W_A))\ra \coder_B(T_B(W_B))$$
given by $f_*(\delta_A)=(1\tns f)\delta_A$, moreover if $\delta_A$ is a codifferential
then so is $f_*(A)$. A codifferential $d_A$ on $T_A(W_A)$ is said to be a deformation
of the codifferential $d$ on $T(W)$ if $\epsilon_*(d_A)=d$.

If $d_A$ is a deformation of $d$ with base $A$ then we can express
\begin{equation*}
 d_A=d+\ph
\end{equation*}
where $\ph:T(W)\ra T(W)\tns\m$. The condition for $d_A$ to be a codifferential is
the Maurer-Cartan equation,
\begin{equation*}
D(\ph)+\frac12[\ph,\ph]=0
\end{equation*}
If $\m^2=0$ we say that $A$ is an infinitesimal algebra and a deformation with
base $A$ is called infinitesimal.

A typical example of an infinitesimal base is $\k[t]/(t^2)$, moreover,
the classical notion of an infinitesimal deformation
$$d_t=d+t\ph$$
is precisely an infinitesimal deformation with base $\k[t]/(t^2)$.

A local algebra $A$ is complete if
\begin{equation*}
 A=\invlim_kA/\m^k
\end{equation*}
A complete, local augmented $\k$-algebra will be called formal
and a deformation with a formal base is called a formal deformation.
An infinitesimal base is automatically formal, so every infinitesimal
deformation is a formal deformation.

An example of a formal base is $A=\k[[t]]$ and a deformation of $d$
with base $A$ can be expressed in the form
$$d_t=d+t\psi_1+t^2\psi_2+\dots$$
This is the classical notion of a formal deformation.
It is easy to see that the condition for $d_t$ to be a formal deformation reduces to
\begin{align*}
 D(\psi_{n+1})=-\frac12\sum_{k=1}^{n}[\psi_k,\psi_{n+1-k}]
\end{align*}

An automorphism of $W_A$ over $A$ is an $A$-linear isomorphism
$g_A:W_A\ra W_A$ making the diagram below commute.
\begin{figure}[h!]
 $$\xymatrix{
 W_A \ar[r]^{g_A} \ar[d]^{\epsilon_*} & W_A \ar[d]^{\epsilon_*} \\
 W \ar[r]^I & W}$$
\end{figure}
The map $g_A$ is induced by its restriction to $T(W)\tns\k$ so we can view $g_A$ as a map
$$g_A:T(W)\ra T(W)\tns A$$
so we ca express $g_A$ in the form
$$g_A=I+\lambda$$
where $\lambda:T(W)\ra T(W)\tns\m$. If $A$ is infinitesimal then $g_A^{-1}=I-\lambda$.

Two deformations $d_A$ and $d_A'$ are said to be equivalent over $A$
if there is an automorphism $g_A$ of $W_A$ over $A$ such that $g_A^*(d_A)=d_A'$.
In this case we write $d'_A\sim d_A$.

An infinitesimal deformation $d_A$ with base $A$ is called universal
if whenever $d_B$ is an infinitesimal deformation with base $B$, there is a unique
morphism $f:A\ra B$ such that $f_*(d_A)\sim d_B$.

\begin{thm}
 If $\dim H^2_{odd}(d)<\infty$ then there is a universal infinitesimal deformation
 $\dinfty$ of $d$. Given by
 $$\dinfty=d+\delta^it_i$$
 where $H^2_{odd}(d)=\langle\bar{\delta^i}\rangle$ and $A=\k[t_i]/(t_it_j)$ is the
 base of deformation.
\end{thm}

A formal deformation $d_A$ with base $A$ is called versal if given any formal
deformation of $d_B$ with base $B$ there is a morphism $f:A\ra B$ such that $f_*(d_A)\sim d_B$.
Notice that the difference between the versal and the universal property
of infinitesimal deformations is that $f$ need not be unique.
A versal deformation is called \emph{miniversal} if $f$ is unique whenever
$B$ is infinitesimal. The basic result about versal deformation is:
\begin{thm}
 If $\dim H^2_{odd}(d)<\infty$ then a miniversal deformation of $d$ exists.
\end{thm}

In this paper we will only need the following result to compute the versal deformations.
\begin{thm}
 Suppose $H^2_{odd}(d)=\langle\bar{\delta^i}\rangle$ and $[\delta^i,\delta^j]=0$
 for all $i,j$ then the infinitesimal deformation
 $$\dinfty=d+\delta^it_i$$
 is miniversal, with base $A=\k[[t_i]].$
\end{thm}

The construction of the moduli space as a geometric object is based on the idea
that codifferentials
which can be obtained by deformations with small parameters are ``close'' to each other.
From the small deformations,
we can construct 1-parameter families or even multi-parameter families, which are defined
for small values of the
parameters, except possibly when the parameters vanish.

If $d_t$ is a one parameter family of deformations, then two things can occur.
First, it may happen that
$d_t$ is equivalent to a certain codifferential $d'$ for every small value of $t$ except zero.
Then we say that $d_t$ is a jump deformation from $d$ to $d'$. It will
never occur that $d'$ is equivalent to $d$,
so  there are no jump deformations from a codifferential to itself.
Otherwise, the codifferentials $d_t$ will all be nonequivalent
if $t$ is small enough. In this case, we
say that $d_t$ is a smooth deformation.

In \cite{fp10}, it was proved for Lie algebras that given three
codifferentials $d$, $d'$ and $d''$,
if there are jump deformations from $d$ to $d'$ and from $d'$ to $d''$,
then  there is a jump deformation from
$d$ to $d''$. The proof of the corresponding statement for associative
algebras is essentially the same.

Similarly, if there is a jump deformation from $d$ to $d'$, and a family of smooth deformations
$d'_t$, then there is a family $d_t$ of smooth deformations of $d$,
such that every deformation in the image of
$d'_t$ lies in the image of $d_t$, for sufficiently small values of $t$.
In this case, we say that the
smooth deformation of $d$ factors through the jump deformation to $d'$.

In the examples of complex moduli spaces of Lie and associative algebras which we have studied,
it turns out that there is a natural
stratification of the moduli space of $n$-dimensional
algebras by orbifolds, where the codifferentials
on a given strata are connected by smooth deformations,
which don't factor through jump deformations.
These smooth deformations determine the local neighborhood structure.

The strata are connected by jump deformations, in the sense that
any smooth deformation from a codifferential on one strata
to another strata factors through a jump deformation.
Moreover, all of the strata are given by projective
orbifolds. In fact, in all the complex examples we have studied,
the orbifolds either are single points,
or $\C\P^n$ quotiented out by either $\Sigma_{n+1}$ or a subgroup, acting on
 $\C\P^n$ by permuting the coordinates.

We don't have a concrete proof at this time, but we conjecture
that this pattern holds in general. In other words,
we believe the following conjecture.
\begin{con}[Fialowski-Penkava]
The moduli space of Lie or associative algebras of a fixed finite
dimension $n$ are stratified by projective orbifolds,
with jump deformations and smooth deformations factoring through
jump deformations providing the only
deformations between the strata.
\end{con}
\section{Deformations of the elements in the moduli space}
\subsection{$d_1=\psa{23}2-\psa{32}2+\psa{22}3-\psa{33}3$}
$$\left[ \begin {array}{ccccccccc}
0&0&0&0&0&0&0&0&0\\\noalign{\medskip}
0&0&0&0&0&1&0&-1&0\\\noalign{\medskip}
0&0&0&1&0&0&0&0&-1
\end {array} \right],$$
This is the direct sum of the $1|1$-dimensional complex simple algebra
 and the trivial $1|0$-dimensional
algebra. This algebra is not unital. Its center is
spanned by $\{v_1, v_3\}$. We have $H^n=\langle\pha{1^n}1\rangle$ for $n>0$.
Since $H^2$ has no odd elements, the algebra is rigid.

\subsection{$d_2=\psa{33}3+\psa{31}1+\psa{32}2$}
The matrix of this codifferential is
$$\left[ \begin {array}{ccccccccc}
0&0&0&0&0&0&1&0&0\\\noalign{\medskip}
0&0&0&0&0&0&0&1&0\\\noalign{\medskip}
0&0&0&0&0&0&0&0&1
\end {array} \right],$$
This algebra is the first of a family of rigid extensions of the $0|1$-dimensional simple
algebras. We have $H^1-\langle\pha22,\pha21,\pha12\rangle$, and $H^n=0$ otherwise.
Its opposite algebra is $d_3$.

\subsection{$d_3=\psa{33}3-\psa{13}1-\psa{23}2$}
The matrix of this codifferential is
$$\left[ \begin {array}{ccccccccc}
0&0&0&0&-1&0&0&0&0\\\noalign{\medskip}
0&0&0&0&0&-1&0&0&0\\\noalign{\medskip}
0&0&0&0&0&0&0&0&1
\end {array} \right],$$
This algebra is the second of a family of rigid extensions of the $0|1$-dimensional simple
algebras. We have $H^1=\langle\pha22,\pha21,\pha12\rangle$, and $H^n=0$ otherwise.
Its opposite algebra is $d_3$.
\subsection{$d_4=\psa{33}3+\psa{31}1-\psa{23}2$}
The matrix of this codifferential is
$$\left[ \begin {array}{ccccccccc}
0&0&0&0&0&0&1&0&0\\\noalign{\medskip}
0&0&0&0&0&-1&0&0&0\\\noalign{\medskip}
0&0&0&0&0&0&0&0&1
\end {array} \right],$$
This algebra is the third of a family of rigid extensions of the $0|1$-dimensional simple
algebras. We have $h^n=1|0$ if $n$ is odd and $h^n=0|0$ otherwise.
The algebra is isomorphic to its opposite algebra.

\subsection{$d_5=\psa{33}3-\psa{13}1+\psa{32}2-\psa{23}2$}
The matrix of this codifferential is
$$\left[ \begin {array}{ccccccccc}
0&0&0&0&0&0&1&0&0\\\noalign{\medskip}
0&0&0&0&0&-1&0&1&0\\\noalign{\medskip}
0&0&0&0&0&0&0&0&1
\end {array} \right],$$
This algebra is the fourth of a family of rigid extensions of the $0|1$-dimensional simple
algebras. We have $H^n=\langle\pha{2^n}2\rangle$ for all $n$.
Its opposite algebra is $d_6$.
\subsection{$d_6=\psa{33}3+\psa{31}1+\psa{32}2-\psa{23}2$}
The matrix of this codifferential is
$$\left[ \begin {array}{ccccccccc}
0&-1&1&0&0&0&0&0&0\\\noalign{\medskip}
0&0&0&1&0&0&0&0&0\\\noalign{\medskip}
0&0&0&0&0&0&0&0&1
\end {array} \right],$$
This algebra is the fifth of a family of rigid extensions of the $0|1$-dimensional simple
algebras. We have $H^n=\langle\pha{2^n}2\rangle$ for all $n$.
Its opposite algebra is $d_5$.
\subsection{$d_7=\psa{33}3+\psa{32}2$}
The matrix of this codifferential is
$$\left[ \begin {array}{ccccccccc}
0&0&0&0&0&0&0&0&0\\\noalign{\medskip}
0&0&0&0&0&0&0&1&0\\\noalign{\medskip}
0&0&0&0&0&0&0&0&1
\end {array} \right],$$
This algebra is the fifth of a family of rigid extensions of the $0|1$-dimensional simple
algebras.
We have
\begin{align*}
H^0&=\langle\pha{}1\rangle\\
H^1&=\langle\pha{1}1\rangle\\
H^{n}&=\langle\pha{1^n}1,\pha{21^{n-1}}2\rangle, \text{ if $n>1$}.
\end{align*}
Since $h^2=2|0$, there are no odd elements in $H^2$, so this algebra is rigid.
Its opposite algebra is $d_8$.
\subsection{$d_8==\psa{33}3-\psa{23}2$}
The matrix of this codifferential is
$$\left[ \begin {array}{ccccccccc}
0&0&0&0&0&0&0&0&0\\\noalign{\medskip}
0&0&0&0&0&-1&0&0&0\\\noalign{\medskip}
0&0&0&0&0&0&0&0&1
\end {array} \right].$$
This algebra is the sixth of a family of rigid extensions of the $0|1$-dimensional simple
algebras.
We have
\begin{align*}
H^0&=\langle\pha{}1\rangle\\
H^1&=\langle\pha{1}1\rangle\\
H^{n}&=\langle\pha{1^n}1,\pha{1^{n-1}2}2\rangle, \text{ if $n>1$}.
\end{align*}
Since $h^2=2|0$, there are no odd elements in $H^2$, so this algebra is rigid.
Its opposite algebra is $d_7$.
\subsection{$d_9=\psa{33}3+\psa{31}1-\psa{13}1+\psa{32}2-\psa{23}2$}
The matrix of this codifferential is
$$\left[ \begin {array}{ccccccccc}
0&0&0&0&-1&0&1&0&0\\\noalign{\medskip}
0&0&0&0&0&-1&0&1&0\\\noalign{\medskip}
0&0&0&0&0&0&0&0&1
\end {array} \right].$$
This algebra is the seventh of a family of rigid extensions of the $0|1$-dimensional simple
algebras. It is both unital and commutative.
Since $h^2=6|0$, there are no odd elements in $H^2$, so this algebra is rigid.
\subsection{$d_{10}=\psa{33}3$}
The matrix of this codifferential is
$$\left[ \begin {array}{ccccccccc}
0&0&0&0&0&0&0&0&0\\\noalign{\medskip}
0&0&0&0&0&0&0&0&0\\\noalign{\medskip}
0&0&0&0&0&0&0&0&1
\end {array} \right].$$
This algebra is the direct sum of the $0|1$-dimensional simple algebra
$\C$ and the trivial nilpotent $2|0$-dimensional algebra.
The algebra is not unital but is commutative. Since $h^2=8|0$, the algebra is rigid.
In fact, except for $H^0$, the cohomology is the subspace of all cochains in $C(M)$. so
$h^n=2^{n+1}|0$ for $n>0$.
\subsection{$d_{11}=\psa{33}3+\psa{32}2-\psa{23}2$}
The matrix of this codifferential is
$$\left[ \begin {array}{ccccccccc}
0&0&0&0&0&0&0&0&0\\\noalign{\medskip}
0&0&0&0&0&-1&0&1&0\\\noalign{\medskip}
0&0&0&0&0&0&0&0&1
\end {array} \right].$$
This algebra is a  extension of the $0|1$-dimensional simple algebra
$\C$ by the trivial $2|0$-dimensional algebra.
The algebra is not unital but is commutative. Its opposite algebra is $d_{10}$.
We have $H^n=\langle\pha{1^n}1,\pha{2^n}2,\psa{2^n}3\rangle$ for all $n$.
The versal deformation is given by $\dinfty=d+\psa{22}3t$, which is a jump
deformation to $d_1$.

\subsection{$d_{12}=\psa{22}3+\psa{23}1-\psa{32}1$}
The matrix of this codifferential is
$$\left[ \begin {array}{ccccccccc}
0&0&0&0&0&0&1&0&0\\\noalign{\medskip}
0&0&0&0&0&0&0&1&0\\\noalign{\medskip}
0&0&0&0&0&0&0&0&1
\end {array} \right].$$
This algebra is an extension of the $1|1$-dimensional nilpotent algebra
$\delta=\psa{22}3$ by the trivial $1|0$-dimensional algebra $\mu=0$.
As is true for all nilpotent algebras, it is not unital. Its center is spanned
by $\{v_1,v_3\}$. We have $h^n=2|0$ when $n$ is odd and $h^n=1|1$ when $n$ is even.
The versal deformation is given by
\begin{equation*}
\dinfty=d+\psa{31}1t-\psa{13}1t+\psa{11}3t^2-\psa{12}3t=\psa{21}3t+\psa{33}3t,
\end{equation*}
which is a jump deformation to $d_1$.
\subsection{$d_{13}=\psa{22}3+\psa{21}3p+\psa{12}3q$}
The matrix of this codifferential is
$$\left[ \begin {array}{ccccccccc}
0&0&0&0&0&0&0&0&0\\\noalign{\medskip}
0&0&0&0&0&0&0&0&0\\\noalign{\medskip}
0&q&p&1&0&0&0&0&0
\end {array} \right].$$
This is a family of nilpotent extensions of the $1|1$-dimensional nilpotent algebra
$\delta=\psa{22}3$ by the trivial $1|0$-dimensional algebra. The family is projectively
parameterized by $(p:q)\in\C\P^1/\Sigma_2$, where the action of $\Sigma_2$ on $\C\P^1$
is given by permuting the coordinates. Thus $d_{13}(p:q)\sim d_{13}(q:p)$. The center
of this algebra is spanned by $\{v_3\}$.
The special points $(1:0)$, $(1:1)$ $(1:-1)$ and $(0:0)$ have different cohomology than
the generic pattern.

Generically $h^2=2|1$ and the versal deformation is given by $\dinfty=d_{13}(p+t:q)$,
which is a smooth
deformation along the family. For the special point $(1:1)$, $h^2=2|1$, and the generic
formula for the versal deformation holds, so it is not special in terms of its
deformations. For the special point $(1:0)$, $h^2=4|1$ and the versal
deformation is given by $\dinfty=d+\psa{11}3t$, which is equivalent to
$d_{13}(1+\sqrt{1-4t}:2t)$. For the special point $(1:-1)$, $h^2=3|1$ the versal deformation
is given by the same formula as for $(1:0)$, but this time it is equivalent to
$d_{13}(1+t:-1+t+\sqrt{-t})$.

Finally, the generic point $d_{13}(0:0)$ has a more interesting deformation theory.
We have $h^2=5|4$, and $h^3=10|7$, so it is not surprising that there are relations
on the base of the versal deformation.
Its versal deformation has matrix
\begin{equation*}
\left[ \begin {array}{ccccccccc}
0&0&0&0&0&-t_{{1}}&t_{{1}}t_{{3}}&t_{{1}}&0\\\noalign{\medskip}
0&0&0&0&-t_{{2}}t_{{3}}&-t_{{2}}+t_{{1}}t_{{3}}
&t_{{1}}t_{{4}}+t_{{2}}t_{{3}}-2\,{t_{{3}}}^{2}t_{{1}}&t_{{2}}-2\,
t_{{1}}t_{{3}}&0\\\noalign{\medskip}t_{{4}}&0&t_{{3}}&1&0&0&0&0&t_{{2}
}\end {array} \right].
\end{equation*}
The third order deformation is versal, and there are 9 nontrivial relations
\begin{align*}
t_1t_4-t_2t_3=0,\quad
t_1(2t_3^2-t_4)=0,\quad
t_1t_3=0\\
(2t_1t_3-t_2)(t_4-t_3^2)=0,\quad
t_1(t_4-t_3)=0,\quad
t_1t_4=0\\
t_1t_3=0,\quad
t_1t_3t_4=0,\quad
t_1t_4+t_2t_3-2t_1t_3^2=0.
\end{align*}
The solutions to these relations are
\begin{equation*}
t_3=t_4=0,\quad t_1=t_2=0,
\end{equation*}
which means the base of the versal deformation is given by two planes that intersect
transversally at the origin.
The first solution gives the codifferential
\begin{equation*}
\dinfty=d+\psa{32}1t_1-\psa{23}1t_1+\psa{32}2t_2-\psa{23}2t_2+\psa{33}3t_2,
\end{equation*}
which is equivalent to $d_1$ except on the line $t_2=0$, where it jumps to $d_{12}$.
The second solution gives the codifferential
\begin{equation*}
\dinfty=d+\psa{21}3t_3+\psa{11}3t_4,
\end{equation*}
which gives jump deformations to $d_{13}(p:q)$ for all $(p:q)$ except $(0:0)$. This behaviour
is consistent with the behaviour of the generic point $(0:0)$ in $\C\P^1$, which is dense
in that space.

\subsection{$d_{14}=\psa{21}3-\psa{12}3$}
The matrix of this codifferential is
$$\left[ \begin {array}{ccccccccc}
0&0&0&0&0&0&0&0&0\\\noalign{\medskip}
0&0&0&0&0&0&0&0&0\\\noalign{\medskip}
0&-1&1&0&0&0&0&0&0
\end {array} \right].$$
This algebra is an extension of the $0|1$-dimensional trivial algebra
$\C$ by the trivial $2|0$-dimensional algebra.
The algebra is commutative.
We have $h^2=6|3$ and $h^3=8|4$, so one might expect there to be relations on the
base of the versal deformation, but in this case there aren't any relations; moreover
the infinitesimal deformation is versal. We have
\begin{equation*}
\dinfty=d+\psa{21}3t_1+\psa{22}3t_2+\psa{11}3t_3.
\end{equation*}
This gives jump deformations to $d_{13}(p:q)$ for all values of $(p:q)$ except $(1:1)$ and
$(0:0)$
\subsection{$d_{15}=\psa{23}1p+\psa{32}1q$}
The matrix of this codifferential is
$$\left[ \begin {array}{ccccccccc}
0&0&0&0&0&p&0&q&0\\\noalign{\medskip}
0&0&0&0&0&0&0&0&0\\\noalign{\medskip}
0&0&0&0&0&0&0&0&0
\end {array} \right].$$
This family of algebras are extensions of the trivial $0|1$-dimensional algebra
$\C$ by the trivial $2|0$-dimensional algebra. They can also be considered
as extensions of the trivial $1|0$-dimensional algebra by the trivial
$1|1$-dimensional algebra. The family is parameterized projectively by
$(p:q)\in\C\P^1$. However, in this case, unlike the $d_{13}(p:q)$ case,
there is no action of the group $\Sigma_2$. In other words, $d_{15}(p:q)$ is
not equivalent to $d_{15}(q:p)$ in general.

The algebra is not commutative unless $p=-q$; generically its center is spanned
by $\{v_1\}$. For the special points $(0:1)$, $(1:0)$, $(1:1)$ and $(1:-1)$, the cohomology
and deformation theory is not generic. In the generic case,
we have  $h^2=1|2$ and $h^3=2|1$, The matrix of the versal deformation is
\begin{equation*}
\left[ \begin {smallmatrix}
0&0&0&0&{\frac {t_{{1}} \left( p-q-t_{{2}} \right) }{2(q+t_{{2}})}}&p&t_{{1}}&q+t_{{2}}&0
\\\noalign{\medskip}
0&0&0&0&
-{\frac {{t_{{1}}}^{2} \left( p+q+t_{
{2}} \right)  \left( p-q-t_{{2}} \right) }{4p \left( q+t_{{2}} \right)
^{2}}}
&-{\frac {t_{{1}} \left( p+q+t_{{2}} \right) }{2(q+t_{{2}}})}
&
-{\frac {{t_{{1}}}^{2} \left( p+q+t_{{2}} \right)  \left( p-q-t_{
{2}} \right) }{ 4\left( q+t_{{2}} \right) {p}^{2}}}&0&0
\\\noalign{\medskip}0&0&0&0&0&0&0&0&t_{{1}}\end {smallmatrix} \right].
\end{equation*}
However, there is one relation:
\begin{equation*}
t_1^2(p+q+t_2)(p-q-t_2)/p^2=0,
\end{equation*}
which in the generic case has only the solution $t_1=0$, which simplifies the
formula for the versal deformation to $\dinfty=d_{15}(p+t_2:q)$. This means that in
the generic case, the deformations are only along the family.

In the $(0:1)$ case, we have $h^2=2|3$ and $h^3=5|3$. The versal deformation is
given by
\begin{equation*}
\dinfty=d+\psa{31}1t_1+\psa{33}3t_1+\psa{23}1t_2+\psa{32}2t_3+\psa{31}1(t_3+t_1t_2).
\end{equation*}
There are two nontrivial relations on the base of the versal deformation,
\begin{equation*}
t_3(t_1+t_3)=0,\quad t_2(t_1+2t_3+t_1t_2)=0,
\end{equation*}
which have solutions
\begin{align*}
t_1=t_3=0,\quad t_2=t_3=0,\quad t_2=0,t_3=-t_1,\quad t_3=0,t_2=-1,\quad t_3=-t_1, t_2=1.
\end{align*}
The last two solutions are not local, in the sense that the lines they parameterize do
not pass through the origin, and thus are not relevant to the deformation picture.
Thus the base of the versal deformation consists of three lines through the origin.
The first line is just $d_{15}(t_2:1)$, which is a deformation along the family, the
second is a jump deformation to $d_7$ and the third is a jump deformation to $d_5$.

In the $(1:0)$ case, we have $h^2=2|3$ and $h^3=5|3$. Note that this algebra is
the opposite algebra to $d_{15}(0:1)$, so its versal deformation could be given
by the opposite algebra to the versal deformation of $d_{15}(0:1)$. In particular,
we have a similar description of the solutions to the versal deformation, and
we obtain deformations along the family, as well as jump deformations to $d_8$ and
$d_6$, the opposite algebras to $d_7$ and $d_5$.

In the $(1:1)$ case, we have $h^2=2|2$ and $h^3=4|2$. The versal deformation
is given by
\begin{equation*}
\dinfty=d+\psa{31}1t_1+\psa{33}3t_1+\psa{32}1t_2-\psa{23}2\frac{t_1(t_2+2)}{2(t_2+1)}
-\psa{31}1\frac{t_1t_2}{2(t_2+1)}+\psa{13}2\frac{t_1^2t_2(t_2+1)}{4(t_2+1)^2}.
\end{equation*}
Note that $\dinfty$ is expressed as a rational function of the parameters rather
than a polynomial. This means that the versal deformation is not given by a finite
order deformation (in terms of this basis of $H^2$), but rather by a power series.
It is very interesting to note that in the examples which we have studied, it has
always been possible to find a basis in which the expression of the versal deformation
is given by a rational expression, although we do not know if this is true in general.

There is one nontrivial relation on the base,
$t_1^2t_2(t_2+2)/(1+t_2)=0$, which has two local solutions,
$t_1=0$ or $t_2=0$. Thus the base of the versal deformation is given by
two lines through the origin. The first line gives $\dinfty=d_{15}(1:1+t_2)$,
which is a deformation along the family, while the second line is a jump deformation
to $d_4$.

In the $(1:-1)$ case, we have $h^2=4|3$ and $h^3=5|4$. We omit the expression
for the versal deformation, because it is a bit nasty. However, there are
two local solutions to the relations on the base of the versal deformation,
$t_2=0$ or $t_1=t_3=0$, which means the base is given by a plane and a line
intersecting transversally at the origin.
The first solution gives
\begin{equation*}
\dinfty=d-\psa{13}1t_1+\psa{31}1t_2+\psa{33}3t_1+\psa{22}3t_3+\psa{11}3t_1^2t_3
-\psa{12}3t_1t_3-\psa{21}3t_1t_3.
\end{equation*}
When neither $t_1$ nor $t_3$ vanish, the deformation is equivalent to $d_1$,
while on the line $t_3=0$, it jumps to $d_{11}$, and on the line $t_1=0$,
it jumps to $d_{12}$.

The second solution is $\dinfty=d_{15}(1:-1+t_2)$, which is just a deformation
along the family.

Note that for this family, the generic point $d_{15}(0:0)$ is just the zero
codifferential, which obviously has jump deformations to every codifferential
in the moduli space.

\bibliographystyle{amsplain}
\providecommand{\bysame}{\leavevmode\hbox to3em{\hrulefill}\thinspace}
\providecommand{\MR}{\relax\ifhmode\unskip\space\fi MR }
\providecommand{\MRhref}[2]{%
  \href{http://www.ams.org/mathscinet-getitem?mr=#1}{#2}
}
\providecommand{\href}[2]{#2}

\end{document}